\documentclass{article} 
\usepackage[left=0.8in,top=0.8in,right=0.8in,bottom=0.8in,head=0.8in]{geometry} 
 \usepackage{url}

\usepackage{hyperref, color}
\usepackage{amsthm}
\usepackage{epsfig, graphics, graphicx}
\usepackage{latexsym}
\usepackage[parfill]{parskip}
\usepackage[tight]{subfigure}
 \usepackage{amsmath,amssymb,enumerate,comment}
\usepackage{accents}
\usepackage[square,numbers]{natbib}

{\begin{enumerate}\item[(#1)]}%
{\end{enumerate}}

{\begin{enumerate}\item[(#1)] [#2 points] }%
{\end{enumerate}}

\newcommand{\E}{\mathbb{E}}

\newcommand{\R}{\mathbb{R}}
\newcommand{\Hk}{\mathcal{H}_k}

\newcommand{\zeros}{\mathbf{0}}

\newcommand{\half}{\tfrac1{2}}

\title{Rows vs Columns for Linear Systems of Equations -\\ Randomized Kaczmarz or  Coordinate Descent?}

\author{
Aaditya Ramdas\\
Machine Learning Department\\
Carnegie Mellon University\\
\texttt{aramdas@cs.cmu.edu} \\
}

\begin{document}

\maketitle

\begin{abstract}
This paper is about randomized iterative algorithms for solving a linear system of equations $X \beta = y$ in different settings. Recent interest in the topic was reignited when Strohmer and Vershynin (2009) proved the linear convergence rate of a Randomized Kaczmarz (RK) algorithm that works on the rows of $X$ (data points). Following that, Leventhal and Lewis (2010) proved the linear convergence of a Randomized Coordinate Descent (RCD) algorithm that works on the columns of $X$ (features). The aim of this paper is to simplify our understanding of these two algorithms, establish the direct relationships between them (though RK is often compared to Stochastic Gradient Descent), and examine the algorithmic commonalities or tradeoffs involved with working on rows or columns. We also discuss Kernel Ridge Regression and present a Kaczmarz-style algorithm that works on data points and having the advantage of solving the problem without ever storing or forming the Gram matrix, one of the recognized problems encountered when scaling kernelized methods. 
\end{abstract}

\section{Introduction}
Solving linear systems of equations is a classical topic and our interest in this problem is fairly limited. While we do compare two algorithms - Randomized Kaczmarz (RK) (see \citet{StrVer09}) and Randomized Coordinate Descent (RCD) (see \citet{LevLew10}) - to each other, we will not be presently concerned with comparing them to the host of other classical algorithms in the literature, like Conjugate Gradient and  Gradient Descent methods. 
Our primary aim will be to understand the algorithmic similarities and differences involved in working with rows and columns (RK and RCD) of a given input.

Assume we have a known $n \times p$ matrix $X$ (representing $n$ data points with $p$ features each), an unknown $p$-dimensional vector $\beta$ (regression coefficients), and a known $n$-dimensional vector $y$ (observations). The system of equations that one wants to solve is
\begin{equation}\label{eq:syseq}
X \beta = y
\end{equation}

When there exists at least one solution to the above system, we say that it is consistent. 
When there is a unique consistent solution, solving \eqref{eq:syseq} is a special case of the more general problem of minimizing the residual norm (makes sense when there are no consistent solutions) 
\begin{equation}\label{eq:LS}
\min_{\beta \in \mathbf{R}^p} \half \|y-X\beta\|^2 =: L(\beta)
\end{equation}

When there is a unique consistent solution, it is also a special case of the more general problem of finding the minimum norm consistent solution (makes sense when there are infinite solutions)
\begin{equation}\label{eq:minnorm}
\min_{\beta \in \mathbf{R}^p} \|\beta\|^2 \mbox{ \ s.t. \ } y = X\beta
\end{equation}

Sparse versions of the above problems, while interesting, and not in the scope of this work. We will later consider the Ridge Regression extension to Eq.\eqref{eq:LS}.

We represent the $i$-th row ($i=1,...,n$) of $X$ by $X^i$, and the $j$-th column ($j=1,...,p$) by $X_j$. Similarly, the $i$-th observation is $y^i$, and the $j$-th regression coefficient is $\beta_j$. The reason for the linear regression setup using statistical notation of $n,p$ and $X,\beta,y$, is simply author comfort, and the literature sometimes uses $Ax=b$ instead, with an $m \times n$ matrix A.

For the rest of this paper, we refer to a method as being Kaczmarz-like when its updates depend on rows (data points) of $X$, like
$$
\beta_{t+1} := \beta_t + \delta_i X^i
$$

where stepsize $\delta_i$ could also depend on $X^i$.
and we refer to a method as being Coordinate Descent style when its updates are on coordinates of $\beta$ and depend on columns (features) of $X$, like
$$
\beta_{t+1} := \beta_t + \delta_j e_j
$$

where stepsize $\delta_j$ could depend on $X_j$. 

Randomized Kaczmarz (and its variations) have been likened to Stochastic Gradient Descent (see \citet{sgdkacz}). Indeed, even before having mentioned the form of $\delta_i$, the update already looks a lot like that of the Perceptron algorithm by \citet{perceptron}, a well known stochastic gradient descent algorithm. However, even though we will later describe some differences from RCD, we argue that RK-style methods are still much more like randomized coordinate descent than stochastic gradient descent algorithms - this is useful not only for interpretation but also for new derivations. We will bring out the intuitive similarity between RK and RCD by establishing striking parallels in proofs of convergence (these proofs are traditionally presented in a different manner), and exploit this in designing a simple RK algorithm for Kernel Ridge Regression. 

There has been a lot of interest and extensions on both these interesting algorithms, and perhaps connections between the two have been floating around in a subtle manner. One of our aims will be to make these connections explicit, intuitively clear, and enable the reader to build an understanding of ideas and proof techniques by the end.  We first need to introduce the aforementioned algorithms before we summarize our contributions. In doing so, we will only refer to papers that are directly relevant to our work.

\subsection{Randomized Kaczmarz (RK)}
Taking  $X, y$ as input and starting from an arbitrary  $\beta_0$, it repeats the following in each iteration. First, pick a random row $r \in \{1...n\}$ with probability proportional to its  Euclidean norm, i.e.
$$
\Pr(r = i) = \frac{\|X^i\|^2}{\|X\|_F^2}
$$
Then, project the current iterate onto that row, i.e.
\begin{equation}
\beta_{t+1} := \beta_t + \frac{(y^r - X^{rT}\beta_t)}{\|X^r\|^2} X^r
\end{equation}

Intuitively, this update can be seen as greedily satisfying the $r$th equation in the linear system, because it is easy to see that after the update,
\begin{equation}
X^{rT}\beta_{t+1} = y^r
\end{equation}

Alternatively, referring to Eq.\eqref{eq:LS}, since  
$
L(\beta) = \half \|y-X\beta\|^2 = \half \sum_{i=1}^n (y^i - X^{iT} \beta)^2,
$ we can interpret this update as stochastic gradient descent (we pick a random data-point on which to update), where the stepsize is the inverse Lipschitz constant of the stochastic gradient  
$$
\nabla^2 \half (y^i - X^{iT} \beta)^2 = \|X^i\|^2.
$$

\citet{StrVer09} showed that the above algorithm 
has an expected linear convergence. We will formally discuss the convergence properties of this algorithm in future sections.

\subsection{Randomized Coordinate Descent (RCD)}
Takeing $X,y$ as input, starting from an arbitrary $\beta_0$, it repeats the following in each iteration. First, pick a random column $c \in \{1...p\}$ with probability proportional to its  Euclidean norm, i.e.
$$
\Pr(c = j) = \frac{\|X_j\|^2}{\|X\|_F^2}
$$

We then minimize the objective $L(\beta) = \half \|y-X\beta\|^2$ with respect to this coordinate to get
\begin{equation}
\beta_{t+1} := \beta_t + \frac{X_c^T(y-X\beta_t)}{\|X_c\|^2} e_c
\end{equation}

where $e_{c}$ is the $c$th coordinate axis. It can be seen as greedily minimizing the objective with respect to the $c$-th coordinate. Indeed, letting $X_{-c},\beta_{-c}$ represent $X$ without its $c$-th column and $\beta$ without its $c$-th coordinate,
\begin{equation}
\frac{\partial L}{\partial \beta_c} = -X_c^T(y-X\beta) = -X_c^T(y-X_{-c}\beta_{-c} - X_c\beta_c)
\end{equation}

Setting this equal to zero for the coordinatewise minimization, we get the aforementioned update for $\beta_c$. Alternately, since $[\nabla L(\beta)]_c = -X_c^T(y-X\beta)$, the above update can intuitively be seen as a univariate descent step where the stepsize is the inverse Lipschitz constant of the gradient along the $c$-th coordinate, since $$[\nabla^2 L(\beta)]_{c,c} = (X^TX)_{c,c} = \|X_c\|^2.$$

\citet{LevLew10} showed that this algorithm has an expected linear convergence. We will discuss the convergence properties of this algorithm in detail in future sections.

\section{Main Results}

We first examine the differences in behavior of the two algorithms in three distinct but related settings. This will bring out the opposite behaviors of the two similar algorithms.

When the system of equations \eqref{eq:syseq} has a unique solution, we represent this by $\beta^*$. This happens when $n \geq p$, and the system is (luckily) consistent. Assuming  that $X$ has full column rank, 
\begin{equation}
\beta^* = (X^TX)^{-1}X^T y
\end{equation}

When \eqref{eq:syseq} does not have any consistent solution, we refer to the least-squares solution of Eq. \eqref{eq:LS} as $\beta_{LS}$. This could happen in the overconstrained case, when $n > p$. Again, assuming  that $X$ has full column rank, we have
\begin{equation}
\beta_{LS} = (X^TX)^{-1}X^T y
\end{equation}

When \eqref{eq:syseq} has infinite solutions, we call the minimum norm solution to \eqref{eq:minnorm} as $\beta^*_{MN}$. This could happen in the underconstrained case, when $n < p$. Assuming  that $X$ has full row rank, we have
\begin{equation}
\beta_{MN}^* = X^T (XX^T)^{-1} y
\end{equation}

In the above notation, the $*$ is used to denote the fact that it is consistent, i.e. it solves the system of equations, the $LS$ stands for Least Squares and $MN$ for minimum norm. We shall return to each of these three situations in that order in future sections of this paper.



One of our main contributions is to achieve a unified understanding of the behaviour of RK and RCD in  these different situations. 
The literature for RK deals with only the first two settings (see  \citet{StrVer09}, \citet{Nee10}, \citet{ZouFre13}), but avoids the third. 
The literature for RCD typically focuses on more general setups than our specific quadratic least squares loss function $L(\beta)$ (see \citet{Nesterov12} or \citet{RicTak12}). However, for both the purposes of completeness, and for a more thorough understanding the relationship between RK and RCD, it turns out to be crucial to analyse all three settings (for equations \eqref{eq:syseq}-\eqref{eq:minnorm}).

\begin{enumerate}
\item When $\beta^*$ is a unique consistent solution, we present proofs of the linear convergence of both algorithms - the results are known from papers by \citet{StrVer09} and \citet{LevLew10} but are presented in  a novel manner so that their relationship becomes clearer and direct comparison is easily possible. 

\item When $\beta_{LS}$ is the (inconsistent) least squares solution, we show why RCD iterates converge linearly to $\beta_{LS}$, but RK iterates do not - making RCD preferable. These facts are not hard to see, but we make it more intuitively and mathematically clear why this should be the case.
\item When $\beta^*_{MN}$ is the minimum norm consistent solution, we explain why RK converges linearly to it, but RCD iterates do not (both so far undocumented observations) - making RK preferable. 
\end{enumerate}

Together, the above three points complete the picture (with solid accompanying intuition) of the opposing behavior of RK and RCD. We then use the insights thus gained to develop a Kaczmarz style algorithm for Ridge Regression and its kernelized version. It is well known that the solution to
\begin{equation}
\min_{\beta \in \R^p} \|y-X\beta\|^2 + \lambda \|\beta\|^2
\end{equation}

can be given in two equivalent forms (using the covariance and gram matrices) as
\begin{equation}
\beta_{RR} ~=~ (X^T X + \lambda I)^{-1}X^T y ~=~ X^T(XX^T + \lambda I)^{-1}y
\end{equation}

The presented algorithms completely avoid inverting, storing or even forming  $XX^T$ and $X^TX$.
Later, we will show that the following updates take only $O(p)$ computation per iteration like RK (starting with $\delta=0, \alpha=\zeros_n, \beta = \zeros_p, r = y$) and have expected linear convergence:
\begin{eqnarray}
\delta_t &=& \frac{y^i - \beta_t^TX^i - \lambda \alpha^i_t}{\|X^i\|^2 + \lambda}\label{eq:RR1}\\
\alpha^i_{t+1} &=& \alpha^i_t + \delta_t\label{eq:RR2} \\
\beta_{t+1} &=& \beta_t + \delta_t X^i\label{eq:RR3}
\end{eqnarray}

where the $i$-th row is picked with probability proportional to $\|X^i\|^2 + \lambda$. 
If $\Hk$ is a Reproducing Kernel Hilbert Space (RKHS, see \citet{learningkernels} for an introduction) associated to positive definite  kernel $k$ and feature map $\phi_x$, it is well known that the solution to the corresponding Kernel Ridge Regression (see \citet{krr}) problem  is 
\begin{eqnarray}
\label{eq:KRRf}
f_{KRR} &=& \arg\min_{f \in \Hk} \sum\limits_{i=1}^n (y_i - f(x_i))^2 + \lambda \|f\|^2_{\Hk} \\
&=& \Phi^T (K+\lambda I)^{-1}y
\end{eqnarray}

where $\Phi = (\phi_{x_1},...,\phi_{x_n})^T$  and $K$ is the gram matrix with $K_{ij}=k(x_i,x_j)$.

One of the main problems with kernel methods is as data size grows, the gram matrix becomes too large to store. This has motivated the study of approximation techniques for such kernel matrices, but we have an alternate suggestion. The aim of a Kaczmarz style algorithm would be to solve the problem by never forming $K$ as exmplified in updates for $\beta_{RR}$. 
For KRR, the update is
\begin{eqnarray}
\alpha^i_{t+1} &=&  \frac{y - \sum_{j\neq i} K(x_i,x_j) \alpha_t^j}{K(x_i,x_i) + \lambda}
\end{eqnarray}

and costs $O(n)$ per iteration, and results in linear convergence as described later. Note that here RK for Kernel Ridge Regression costs $O(n)$ per iteration and RK for Ridge Regression cost $O(p)$ per iteration due to different parameterization. In the latter, we can keep track of $\beta_t$ as well as $\alpha_t$ easily, see Eq.\eqref{eq:RR2},\eqref{eq:RR3}, but for KRR, calculations can only be performed via evaluations of the kernel only ($\beta_t$ corresponds to a function and cannot be stored), and hence have a different cost.

The aforementioned updates and their convergence can  be easily derived after we develop a clear understanding of how RK and RCD methods relate to each other and jointly to positive semi-definite systems of equations. We shall see more of this in Sec.\ref{sec:RR}.

\newpage

\section{Overconstrained System, Consistent}

To be clear, here we will assume that $n > p$, $X$ has full column rank, and that the system is consistent, so
 $y = X\beta^*$.
First, let us write the updates used by both algorithms in a revealing fashion. If RK and RCD picked row $i$ and column $j$ at step $t+1$, and $e^i$ is $1$ in the $i$-th position and $0$ elsewhere, then the updates can be rewritten as below:
\begin{align}\label{eq:RKupdate}
&\mbox{(RK)}& \beta_{t+1} &:= \beta_t + \frac{e^{iT} r_t}{\|X^i\|^2} X^i&\\
&\mbox{(RCD)}& \beta_{t+1} &:= \beta_t + \frac{X_j^T r_t}{\|X_j\|^2}e_j &
\end{align}
where $r_t = y - X\beta_t = X\beta^*-X\beta_t$ is the residual vector. Then multiplying both equations by $X$ gives
\begin{align}
&\mbox{(RK)}&X\beta_{t+1} &:= X\beta_t + \frac{X^{iT}(\beta^*-\beta_t)}{\|X^i\|^2} X X^i&\\
&\mbox{(RCD)}& X\beta_{t+1} &:= X\beta_t + \frac{X_j^T X(\beta^*-\beta_t)}{\|X_j\|^2}X_j& \label{eq:RCDupdate}
\end{align}
We now come to an important difference, which is the key update equation for RK and RCD. 

Firstly, from the update Eq.(\ref{eq:RKupdate}) for RK, we have $\beta_{t+1}-\beta_t$ is parallel to $X^i$. Also, $\beta_{t+1}-\beta^*$ is orthogonal to $X^i$ (why? because $X^{iT}(\beta_{t+1}-\beta^*)=y^i - y^i = 0$). By Pythagoras,
\begin{equation}\label{eq:RKrecursion}
\|\beta_{t+1} - \beta^*\|^2 = \|\beta_t - \beta^*\|^2 - \|\beta_{t+1} - \beta_t\|^2
\end{equation}
Note that from the update Eq.(\ref{eq:RCDupdate}), we have $X\beta_{t+1} - X\beta_t$ is parallel to $X_j$. Also, $X\beta_{t+1} - X\beta^*$ is orthogonal to $X_j$ (why? because $X_j^T(X\beta_{t+1} - X\beta^*) = X_j^T(X\beta_{t+1} - y) = 0$ by the optimality condition $\partial L/\partial \beta_{j} = 0$). By Pythagoras,
\begin{equation}\label{eq:RCDrecursion}
\|X\beta_{t+1} - X\beta^*\|^2 = \|X\beta_t - X\beta^*\|^2 - \|X\beta_{t+1} - X\beta_t\|^2
\end{equation}
The rest of the proof follows by simply substituting for the last term in the above two equations, and is presented in the following table for easy comparison. Note $\Sigma=X^TX$ is the full-rank covariance matrix and we first
 take expectations with respect to the randomness at the $t+1$-st step, conditioning on all randomness up to the $t$-th step. We later iterate this expectation.
\begin{table}[h]
\begin{tabular}{|l|l|}
\hline
Randomized Kaczmarz $\E\|\beta_{t+1} - \beta^*\|^2$ & Randomized Coordinate Descent $\E\|X\beta_{t+1} - X\beta^*\|^2$ \\
\hline 
$= \|\beta_t - \beta^*\|^2 - \E\|\beta_{t+1} - \beta_t\|^2$ & $= \|X\beta_t - X\beta^*\|^2 - \E\|X\beta_{t+1} - X\beta_t\|^2$\\
$= \| \beta_{t} - \beta^*\|^2 -  \sum_i \frac{\|X^i\|^2}{\|X\|_F^2} \frac{(X^i(\beta_t - \beta^*))^2}{(\|X^i\|^2)^2} \|X^i\|^2$ & $= \|X \beta_{t} - X\beta^*\|^2  - \sum_j \frac{\|X_j\|^2}{\|X\|_F^2} \frac{(X_j^T X(\beta_t - \beta^*))^2}{(\|X_j\|^2)^2}\|X_j\|^2$ \\
$= \|\beta_{t} - \beta^*\|^2 \left(1 - \frac1{\|X\|_F^2} \frac{\| X(\beta_t - \beta^*)\|^2}{\| \beta_{t} - \beta^*\|^2} \right)$ & $= \|X \beta_{t} - X\beta^*\|^2 \left(1 - \frac1{\|X\|_F^2}  \frac{\|X^T X(\beta_t - \beta^*)\|^2}{\|X \beta_{t} - X\beta^*\|^2} \right)$ \\
$\leq \|\beta_t-\beta^*\|^2 (1 - \frac{\sigma_{\min}(\Sigma)}{Tr(\Sigma)})$ & $\leq \|X\beta_t-X\beta^*\|^2 (1 - \frac{\sigma_{\min}(\Sigma)}{Tr(\Sigma)})$ \\
\hline
\end{tabular}
\end{table}

Here, $\sigma_{\min}(\Sigma)\| \beta_t - \beta^*\|^2 \leq \|X(\beta_t - \beta^*)\|^2 $  i.e. $\sigma_{\min}(\Sigma)$ is the smallest eigenvalue. 
It follows that
\begin{align}
&\mbox{(RK)}& \E\|\beta_t - \beta^*\|^2 &\leq \left( 1 - \frac{\sigma_{\min}(\Sigma)}{Tr(\Sigma)} \right)^{t}\|\beta_0-\beta^*\|^2&\\ \label{eq:RCDlin}
&\mbox{(RCD)}& \E\|\beta_t - \beta^*\|_\Sigma^2 &\leq \left( 1 - \frac{\sigma_{\min}(\Sigma)}{Tr(\Sigma)} \right)^{t}\|\beta_0-\beta^*\|_\Sigma^2 
\end{align}

Since $\Sigma$ is invertible when $n>p$ and $X$ has full column rank, the last equation also implies linear convergence of $\E\|\beta_t-\beta^*\|^2$.
The final results do exist in \citet{StrVer09,LevLew10} but there is utility in seeing the two proofs in a form that differs from their original presentation, side by side. In this setting, both RK and RCD are essentially equivalent (without computational considerations).

\newpage

\section{Overconstrained System, Inconsistent}

Here, we will assume that $n > p$, $X$ is full column rank, and the system is (expectedly) inconsistent, so $y = X\beta_{LS} + z$, where $z$ is such that $X^T z = 0$. It is easy to see this condition, because as mentioned earlier,
\[
\beta_{LS} = (X^T X)^{-1}X^T y
\]
implying that $X^T X \beta_{LS} = X^T y$. Substituting $y=X\beta_{LS}+z$ gives that $X^Tz=0$.

In this setting, RK is known to not converge to the least squares solution, as is easily verified experimentally. The tightest convergence upper bounds known are by \citet{Nee10} and \citet{ZouFre13} who show that
\[
\E\|\beta_t - \beta_{LS}\|^2 \leq \left( 1 - \frac{\sigma_{\min}^+(\Sigma)}{Tr(\Sigma)} \right)^{t}\|\beta_0-\beta_{LS}\|^2 + \frac{\|w\|^2}{\sigma_{\min}^+(\Sigma)^2}
\]
If you tried to follow the previous proof, Eq.\eqref{eq:RKrecursion} does not go through - Pythagoras fails because $\beta_{t+1}-\beta_{LS}\not\perp X^i$, since $X^{iT}(\beta_{t+1}-\beta_{LS}) = y^i-X^{iT}\beta_{LS} \neq 0 $.
Intuitively, the reason RK does not converge  is that every update of RK (say of row $i$) is a projection onto the ``wrong'' hyperplane that has constant $y^i$ (where the ``right'' hyperplane would involve projecting onto a parallel hyperplane with constant $y^i-z^i$ where $z$ was defined above). An alternate intuition is that all RK updates are in the span of the rows, but $\beta_{LS}$ is not in the row span. These intuitive explanations are easily confirmed by experiments seen in \citet{Nee10,ZouFre13}.

The same doesn't hold for RCD. Almost magically, in the previous proof, Pythagoras works because 
$$X_j^T(X\beta_{t+1} - X\beta_{LS}) = X_j^T(X\beta_{t+1} - y) + X_j^T(y-X\beta_{LS}) = 0$$
The first term is 0 by the  optimality condition for $\beta_{t+1}$ i.e. $X_j^T(X\beta_{t+1} - y) = \partial L/\partial \beta_j = 0$. The second term is zero by the global optimality of $\beta_{LS}$ i.e. $X^T(y - X\beta_{LS}) = \nabla L = 0$. Also, $\Sigma$ is full rank as before.
Hence, since RCD works in the space of fitted values $X\beta$ and not the iterates $\beta$.

In summary, RK does not converge to the LS solution, but RCD does at the same linear rate.
The Randomized Extended Kaczmarz by \citet{ZouFre13} is a modification of RK designed to converge by randomly projecting out $z$, but its discussion is beyond our  scope. We were also alerted to a recent, independent Arxiv paper by \citet{frek}.

\section{Underconstrained System, Infinite Solutions}

Here, we will assume that $p > n$, $X$ is full row rank and the system is (expectedly) consistent with infinite solutions. As mentioned earlier, it is easy to show that
\[
\beta^*_{MN} = X^T(XX^T)^{-1}y
\]
(which clearly satisfies $X\beta^*_{MN} = y$). Every other consistent solution can be expressed as 
$$
\beta^* = \beta^*_{MN} + z ~\mbox{~ where ~}~ Xz=0
$$ 
Clearly any such $\beta^*$ would also satisfy $X\beta^* = X\beta^*_{MN}=0$. Since $Xz=0$, $z \perp \beta^*_{MN}$ implying $\|\beta^*\|^2 = \|\beta^*_{MN}\|^2 + \|z\|^2$, showing that $\beta^*_{MN}$ is indeed the minimum norm solution as claimed.

In this case, RK has good behaviour, and starting from $\beta_0=0$, it does converge linearly to $\beta^*_{MN}$. Intuitvely, $\beta^*_{MN} = X^T \alpha$ (for $\alpha = (XX^T)^{-1}y$) and hence is in the row span of $X$. Starting from $\beta_0=0$, RK only adds multiples of rows to its iterates, and hence will never have any component orthogonal to the row span of $X$ (i.e. will never add any $z$ such that $Xz=0$). There is exactly one solution with no component orthogonal to the row span of $X$, and that is $\beta^*_{MN}$, and hence RK converges linearly to the required point. It is important not to start from an arbitrary $\beta_0$ since the RK updates can never wipe out any component of $\beta_0$ that is perpendicular to the row span of $X$.

Mathematically, the previous earlier proof works because Pythagoras goes through since it is a consistent system.  However, $\Sigma$ is not full rank but note that since both $\beta^*_{MN}$ and $\beta_t$ are in the row span, $\beta_t - \beta^*_{MN}$ has no component orthogonal to $X$ (unless it equals zero and we're done). Hence $\sigma_{\min}^+(\Sigma)\| \beta_t - \beta^*\|^2 \leq \|X(\beta_t - \beta^*)\|^2$ does hold $\sigma_{\min}^+$ being the smallest positive eigenvalue of $\Sigma$.

RCD unfortunately suffers the opposite fate. The iterates do not converge to $\beta_{LS}^*$, even though $X\beta_t$ does converge to $X\beta^*$. Mathematically, the convergence proof still carries forward as before till Eq.\eqref{eq:RCDlin}, but in the last step where $X^T X$ cannot be inverted because it is not full rank. Hence we get convergence of the residual to zero, without getting convergence of the iterates to the least squares solution.

Unfortunately, when each update is cheaper for RK than RCD (due to matrix size), RCD is preferred for  reasons of convergence and when it is cheaper for RCD than RK, RK is preferred.




\section{Randomized Kaczmarz for Ridge Regression}\label{sec:RR}

Both RK and RCD can be viewed in the following fashion. Suppose we have a positive definite matrix $A$, and we want to solve $Ax=b$. Instead of casting it as $\min_x \|Ax-b\|^2$, we can alternatively pose the different problem $\min_x \half x^TAx - b^Tx$. Then one could use the update 
\[
x_{t+1} = x_t + \frac{b_i-A_i^Tx_t}{A_{ii}} e_i
\]
where the $b_i-A_i^T x_t$ is basically the $i$-th coordinate of the gradient, and $A_{ii}$ is the Lipschitz constant of the $i$-th coordinate of the gradient (related works include \citet{LevLew10}, \citet{Nesterov12},\citet{RicTak12},\citet{SidLee13}).

In this light, the RK update can be seen as the randomized coordinate descent rule for the psd system $XX^T \alpha = y$ (substituting $\beta = X^T \alpha$) and treating $A = XX^T$ and $b=y$.

Similarly, the RCD update can be seen as the randomized coordinate descent rule for  the psd system $X^T X \beta = X^T y$ and treating $A = X^T X$ and $b=X^T y$.

Using this connection, we propose the following update rule:
\begin{eqnarray}
\delta_t &=& \frac{y^i - \beta_t^TX^i - \lambda \alpha^i_t}{\|X^i\|^2 + \lambda} \label{eq:kacz1}\\
\alpha^i_{t+1} &=& \alpha^i_t + \delta_t  \label{eq:kacz2}\\
\beta_{t+1} &=& \beta_t + \delta_t X^i \label{eq:kacz3}
\end{eqnarray}
where the $i$-th row was picked with probability proportional to $\|X^i\|^2 + \lambda$. If all rows are normalized, then this is still a uniform distribution. However, it is more typical to normalize the columns in statistics, and hence one pass over the data must be made to calculate row norms.

We argue that this update rule exhibits linear convergence for $\min_\beta \|y - X\beta\|^2 + \lambda \|\beta\|^2$. Similarly, for Kernel Ridge Regression as mentioned in Eq.\eqref{eq:KRRf}, since one hopes to calculate $f_{KRR} = \Phi^T(K + \lambda I_n)^{-1} y$, the RK-style update can be derived from the randomized coordinate descent update rule for the psd system
\[
(K + \lambda I) \alpha = y
\]
by setting $f_{KRR}=\Phi^T \alpha$. The update for $\alpha$ looks like (where $S_a(z) = \frac{z}{1+a}$)
\begin{eqnarray}\label{eq:kerkacz1}
\alpha^i_{t+1} &=& \frac{K(x_i,x_i)}{K(x_i,x_i)+\lambda}\alpha^i_t + \frac{y_i - \sum_j K(x_i,x_j) \alpha_t^j}{K(x_i,x_i) + \lambda}\\
&=& S_{\frac{\lambda}{K(x_i,x_i)}} \left( \alpha_t^i + \frac{r_i}{K(x_i,x_i)} \right)\label{eq:kerkacz3}
\end{eqnarray}
where row $i$ is picked proportional to $K(x_i,x_i)+\lambda$ (uniform for translation invariant kernels).

Let us contrast this with the randomized coordinate descent update rule   for the loss function $\min_x \half \beta^T (X^TX + \lambda I_p) \beta - y^T X\beta$ i.e. the system $(X^T X + \lambda I_p)\beta = X^T y$.
\begin{eqnarray}
\beta^i_{t+1} &=&\beta^i_t + \frac{X_i^T y - X_i^T X \beta -\lambda \beta_i}{\|X_i\|^2 + \lambda}\label{eq:rcdridge1}\\
&=& \frac{\|X_i\|^2}{\|X_i\|^2 + \lambda}\beta^i_t + \frac{X_i^T r_t}{\|X_i\|^2 + \lambda} \\
&=& S_{\frac{\lambda}{\|X_i\|^2}} \left(\beta^i_t + \frac{X_i^T r_t}{\|X_i\|^2} \right)\label{eq:rcdridge3}
\end{eqnarray}

\subsection{Computation and Convergence}
The RCD updates in \eqref{eq:rcdridge1}-\eqref{eq:rcdridge3} take $O(n)$ time, since each column (feature) is of size $n$. In contrast, the proposed RK updates in \eqref{eq:kacz1}-\eqref{eq:kacz3} takes $O(p)$ time since that is the length of a data point. 

Lastly the RK updates in \eqref{eq:kerkacz1}-\eqref{eq:kerkacz3} take $O(n)$ time (to update $r$) not counting time for kernel evaluations. The difference between the two RK updates for Ridge Regression and Kernel Ridge Regression is that for KRR, we cannot maintain $\alpha$ and $\beta$ since the $\beta$ is a function in the RKHS. This different parameterization makes the updates to $\alpha$ cost $O(n)$ instead of $O(p)$.

While the RK and RCD algorithms are similar and related, one should not be tempted into thinking their convergence rates are the same. Indeed, with no normalization assumption, using a similar style proof as presented earlier, one can show that the convergence rate of the RK for Kernel Ridge Regression is 
\begin{eqnarray}
\E\|\alpha_t - \alpha^*\|_{K+\lambda I_n}^2 &\leq& \left( 1 - \frac{\sigma_{\min}(K + \lambda I_n)}{Tr(K+\lambda I_n)}\right)^{t} \|\alpha_0 - \alpha^*\|_{K+\lambda I_n}^2\\
&=& \begin{cases} \left( 1 - \frac{\lambda}{\sum_i \sigma_i^2 + n\lambda}\right)^{t} \|\alpha_0 - \alpha^*\|_{K+\lambda I_n}^2 ~\mbox{~if $n>p$ }\\
\left( 1 - \frac{\sigma_1^2 + \lambda}{\sum_i \sigma_i^2 + n\lambda}\right)^{t} \|\alpha_0 - \alpha^*\|_{K+\lambda I_n}^2 ~\mbox{~if $p>n$} 
\end{cases}
\end{eqnarray}
and the rate of convergence of the RCD for Ridge Regression is subtly different:
\begin{eqnarray}
\E\|\beta_t - \beta^*\|_{\Sigma+\lambda I_p}^2 &\leq& \left( 1 - \frac{\sigma_{\min}(\Sigma + \lambda I_p)}{Tr(\Sigma + \lambda I_p)}\right)^{t} \|\beta_0 - \beta^*\|_{\Sigma+\lambda I_p}^2 \\
&=& \begin{cases} \left( 1 - \frac{\sigma_1^2 + \lambda}{\sum_i \sigma_i^2 + p\lambda}\right)^{t} \|\beta_0 - \beta^*\|_{\Sigma+\lambda I_p}^2 ~\mbox{~if $n>p$} \\
\left( 1 - \frac{ \lambda}{\sum_i \sigma_i^2 + p\lambda}\right)^{t} \|\beta_0 - \beta^*\|_{\Sigma+\lambda I_p}^2 ~\mbox{~if $p>n$}
\end{cases}
\end{eqnarray}

Kernel Ridge Regression is used as a subroutine in many kernelized machine learning problems (for example, see section 4.3 of \citet{FukSonGre14} for a recent application to a nonparametric state space models). One of the major issues involved with scaling kernel methods is the formation of the gram matrix, which could be prohibitively large to form, store, invert, etc. Our RK-style algorithm gets around this issue completely by never forming the kernel matrix, just as RCD avoids forming $\Sigma$, making it a great choice for scaling kernel methods.



\section{Conclusion}

In this paper, we studied the close connections between the RK and RCD algorithms that have both received a lot of recent attention in the literature, and which we show are in some sense instances of each other when appropriately viewed. While RK is often viewed as a stochastic gradient algorithm, we saw that its ties to RCD are much stronger and that it is easier to understand its convergence through the RCD perspective than the SGD viewpoint.

We first analysed their opposite behavior with linear systems. If the system was consistent and unique then we showed that both algorithms approached the solution at (approximately) the same linear rate, with extremely similar proofs presented for direct comparison. However, if the system was consistent with infinite solutions then we saw that RK converged to the minimum norm solution but RCD didn't, making RK preferable. In contrast, if the system was inconsistent, RCD converged to the least squares solution but RK did not, making RCD the preferred choice. Unfortunately, in both cases, the preferred choices have costlier updates. 

We then exploited the connection between RK and RCD to design an RK-style algorithm for Kernel Ridge Regression (KRR) which avoided explicitly forming and inverting the potentially large kernel matrix. We anticipate this and other randomization techniques will help the scalability of kernel methods.
 

\subsection*{Acknowlegements}
The authors would like to thank 
Deanna Needell and Anna Ma
for some discussions, 
Peter Richtarik   for pointing out \citet{frek}, and
Ryan Tibshirani, for memorable class notes on coordinate descent and for corrections on an earlier version of this manuscript.

\bibliographystyle{agsm}
\bibliography{rvc}

\newpage
\appendix

\end{document}